\documentclass[russian,12pt]{article}
\usepackage{graphicx}
\usepackage{color}
\usepackage{rotating}
\newtheorem{thm}{Theorem}[section]

\newtheorem{al}[thm]{Algorithm}

\begin{document}

\title{\bf New Methods for Solving  Large Scale Linear Programming Problems in the Windows and Linux computer operating systems }

\author{\small Saeed Ketabchi \footnote{Department of Applied Mathematics, Faculty of Mathematical Sciences,
 University of Guilan, P.O. Box 1914, Rasht, Iran,
sketabchi@guilan.ac.ir},\ \
 Hossein Moosaei \footnote{Corresponding author, Department of Mathematics, University of Bojnord, Bojnord, Iran, Email:
 hmoosaei@gmail.com,}  , Hossein Sahleh  \footnote{ Department of Pure Mathematics, Faculty of Mathematical Sciences,  University of Guilan, Email:  sahleh@guilan.ac.irm}~ and Mohammad Hedayati  \footnote{ Department of Pure Mathematics, Faculty of Mathematical Sciences,  University of Guilan, Email:  mb.hedayati@gmail.com}}
\maketitle
\textbf{This paper will be published in Appl. Math. Inf. Sci. Vol. 7 No. 1 (2013)
}
\begin{abstract}
In this study, calculations necessary to solve the large scale linear programming problems in two operating systems, Linux and Windows 7 (Win), are compared using two different methods. Relying on the interior-point methods, linear-programming interior point solvers (LIPSOL) software was used for the first method and relying on an augmented Lagrangian method-based algorithm, the second method used the generalized derivative. The performed calculations for various problems show the produced random in the Linux operating system (OS) and Win OS indicate the efficiency of the performed calculations in the Linux OS in terms of the accuracy and using of the optimum memory.
\end{abstract}
\section{Introduction}
\label{sec.1}
Linear programming (LP) is an important class of optimization problems and is used extensively
in economics, operations research, engineering, and many other fields.

In this paper we consider the primal linear programming  in the standard form

\begin{eqnarray}
 \nonumber
f_*=\min_{x\in X}c^Tx,\quad X=\{x\in R^n:Ax=b,~x\geq 0\},  \qquad  \qquad (P) \\ \nonumber
\end{eqnarray}
where $A\in R^{m\times n},c\in R^n$, and $b\in R^m$ are given, $x$ is primal variable, $0_i$  denotes the i-dimensional zero vector.
In this work we present two  methods for solving problem (P). The first method -which is a Matlab-based package for solving linear programs by interior-Point methods- is LIPSOL (the function linprog in Matlab). The second method  is based on  augmented Lagrangian method. \cite {E0, E2, jin, zhi}

Initially, some problems with different sizes are produced randomly through providing generated problem and then calculations are done in the Win OS. Though the second algorithm enjoys a higher efficiency than MATLAB LINPROG according to the results, it cannot solve the large-scale problems. Now, the operating system is changed from the Win to the Linux and the produced problems are solved using both methods in the Linux OS.
It can be observed that the unsolved problems by the Win OS are solved by the Linux OS. These results show the higher efficiency and flexibility of the Linux OS than the Win OS to make calculations and to solve large scale problems using various methods.

This paper is organized as follows. Augmented Lagrangian method is discussed in the section 2. In Section 3, some examples on various randomly generated problems are provided to illustrate the efficiency and validity of our proposed method. Concluding remarks are given in Section 4.

 In this work by $A^\top$ and $\|.\|$ we mean the transpose of matrix  $A$  and Euclidean norm respectively and $a_+$ replaces negative components of the vector $a$,  by zeros.
\section{Augmented Lagrangian Method }
In this section we consider problem (P).~ Assume that the solution set ~$X_*$ of primal problem ~ (P) is nonempty, hence the solution set $U_*$ of dual problem  is also nonempty, and $x \in R^n$ be an arbitrary vector. Next Theorem tells us that we can get a solution of the dual problem of $(P)$ from the unconstrained minimization problem .
\begin{thm}
\label{QP}
Assume that the solution set $X_*$ of problem (P) is nonempty. Then there exists $\alpha_*>0$ such that for all $\alpha \geq \alpha_*$ the unique least $2-$norm projection $x$ of a point $\bar{x}$ onto $X_*$ is given by
${x}=(\bar{x} +\alpha(A^Tu(\alpha)- c))_+$ where $u(\alpha)$ is a point attaining the minimum in the following problem:
  \begin{eqnarray}
\label{2b}
\min_{u\in R^m}\Phi(u,\alpha,\bar{x}),
\end{eqnarray}
 where $ \Phi(u,\alpha,\bar{x})=-b^Tu+ \frac{1}{2\alpha} \parallel(\bar{x}+\alpha(A^Tu-c))_+ \parallel^2$. In addition,  for all $\alpha> 0$ and $x \in X_*$,  the solution of the  convex,  quadratic problem (\ref{2b}),
 $u_*=u(\alpha) $  is  an exact solution of the dual problem  i.e.  $u(\alpha) \in U_*$.\\
\end{thm}
The proof  is given in \cite{E0}.\\

Now, we describe how can be solved the unconstrained optimization problem (\ref{2b}). The function $\Phi(u,\alpha,\bar{x})$ is augmented Lagrangian function for the  dual of the linear programming $(P)$ (see \cite{A1}),
 \begin{eqnarray}
 \nonumber
f_*=\max_{u\in U}b^Tu,\quad U=\{u\in R^m:A^Tu\leq c\}~~~~~~~~\qquad(D)
\end{eqnarray}

The function $\Phi(u,\alpha,\bar{x})$  is  piecewise quadratic, convex, and just has  the first derivative, but it is not twice  differentiable. Suppose $u$ and $s\in R^m$, for gradient of  $\Phi(u,\alpha,\bar{x})$ we have
    $$  \|\nabla \Phi(u,\alpha,\bar{x}) -\nabla \Phi(s,\alpha,\bar{x})\|\leq \|A\|\|A^T\| \|u-s\|, $$
this means $\nabla \Phi$ is globally Lipschitz continues with constant $K= \|A\|\|A^T\|.$ Thus  for this function generalized Hessian exist and is defined the ~ $m\times m$~ symmetric positive semidefinite matrix \cite{H,K1,M1,M2}
$$\nabla^2\Phi(u,\alpha,\bar{x})=AD(z)A^T,$$
where $D(z)$ denotes an $n\times n$ diagonal matrix with i-diagonal element $z_i$ equals to $1$ if $(\bar{x}+\alpha(A^Tu(\alpha)- c)_i>0$ and equal to $0$ otherwise. Therefore we can use generalized Newton method for solving this problem and to obtain global termination we must use a line-search algorithm (see \cite{N}). In the following algorithm we apply  the generalized Newton method with a line-search based on the Armijo rule \cite{A2}.
\\-----------------------------------------------------------------------------------------------------------------\\
\begin{al}\textbf{Generalized Newton method  with the Armijoo rule}\\
\label{QP1}
Choose any $u_0\in R^m$  and $\epsilon\geq0$\\
 i=0;\\
  \textbf{while}~ \textit{$\|{\nabla \Phi(u_i)}_{\infty}\| \ge   \epsilon$ }\\
 \textbf{Choose} \textit{$\alpha_i = $max$\{s, s\delta, s\delta^2,... \}$ such that} \\
\textit{$\Phi(u_i)-\Phi(u_i+\alpha_id_i)\geq -\alpha_i\mu\nabla \Phi(u_i)d_i, $\\
where $d_i=-\nabla^2 \Phi(u_i)^{-1}\nabla \Phi(u_i)$, $s>0$ be a constant, $\delta\in(0,1) $ and $\mu\in(0,1).$\\
$u_{i+1}=u_i+\alpha_id_i$\\ $i=i+1$;\\}
\textbf{end}
\end{al}
-------------------------------------------------------------------------------------------------------------------\\
In this algorithm, the generalized Hessian may be singular, thus  we used a modified Newton direction as following:
$$-(\nabla^2 \Phi(u_i)+\delta I_m)^{-1}\nabla \Phi(u_i),$$
where $\delta$ is a small positive number $(\delta=10^{-4}),$ and $I_m$ is the identity matrix of m order.\\
Now we  introduce the following iterative process (multiplies method for the dual LP problem (D)):
\begin{eqnarray}
\label{4f}
u^{k+1}=arg\min_{u\in R^n}\{-b^Tu+ \frac{1}{2\alpha} \parallel(x^k+\alpha(A^Tu-c))_+ \parallel^2,
\end{eqnarray}
\begin{eqnarray}
\label{5f}
x^{k+1}=(x^k+\alpha(A^Tu^{k+1}-c))_+,
\end{eqnarray}
where $x^0$ is an arbitrary starting point and the solution  of the problem (\ref{4f}) has been obtained by  algorithm (\ref{QP1}). \\

\begin{thm}
 \label{QP2}
Let the solution set $X_*$ of the problem $(P)$ be nonempty. Then, for all $\alpha>0$ and an arbitrary initial $x^0$ the iterative process (\ref{4f}), (\ref{5f}) converges to $x_*\in X_*$ in finite number of step $k$ and
the primal normal solution $\widehat{x}_*$  was obtained after the first  iteration from above  process, i.e. $k=1$. Furthermore, $u_*=u^{k+1}$ is an exact solution of the dual problem $(D)$.
\end{thm}
The proof of the finite global convergence is given in \cite{A1}

\par\vspace{10pt}
\section{Numerical results}
~~~~In this section  we present  some  numerical results  on various randomly generated problems to the problem (P). The problems are generated using the following MATLAB code:
\\-----------------------------------------------------------------------------------------------------------------\\
\vspace{5pt}
\\ {\tt
\hspace{-.61cm}
\%lpgen: Generate random solvable lp: min c'x s.t. Ax = b ;x>=0,\\
\%Input: m,n,d(ensity); Output: A,b,c; (x, u): primal-dual solution \\
m=input('Enter m:')\\
n=input('Enter n:')\\
d=input('Enter d:')\\
pl=inline('(abs(x)+x)/2')\\
A=sprand(m,n,d);A=100*(A-0.5*spones(A));
x=sparse(10*pl(rand(n,1)));\\
u=spdiags((sign(pl(rand(m,1)-rand(m,1)))),0,m,m)*(rand(m,1)-rand(m,1));\\
b=A*x;c=A'*u+spdiags((ones(n,1)-sign(pl(x))),0,n,n)*10*ones(n,1);\\
format short e;[norm(A*x-b), norm(pl(A'*u-c)), c'*x-b'*u]. \\
}
--------------------------------------------------------------------------------------------------------------------\\
The test problem generator generates a random matrix A for a given $m,n$ and density $d$ and the vector $b$. The elements of A are uniformly distributed between $-50$ and $+50$. In  all  computations were performed,  our variant augmented Lagrangian method  and the generalized  Newton method and Armijo line search, were implemented in MATLAB code. We used  \emph{Core 2 Duo 2.53 GHz} with main memory \emph{4 GB}. The computation results are shown in  Table(1) and  Table(2). We present comparison between the method is based on \emph{\textbf{LIPSOL}} ( see \cite{Z}) in MATLAB  ({\bf linprog})  and our algorithm ({\bf lpf }) in Windows and Linux . \\ Computational results show that, some of the unsolved problems in the Win OS are solved in the Linux OS.
\\
The starting vector used in all following examples is $x^0=0.$ In all solved examples $\alpha=10/{d^{0.5}},~  tol=10^{-10}.$ The total times of computations in each example is given the third   column of the table. The accuracy of optimality conditions of LP problems are in the last three columns.\\

\begin{sidewaystable}
\caption{ Comparison of \textit{\textbf{linprog }} and  \textit{\textbf{lpf}} in  Win}
\begin{center}{
\begin{tabular}{|c|c|c|c|c|c|c|}\hline &&&&\\ [-6mm]
   $m,n,d$&&$ time(sec)$&$\|{x^*}\|$& $\|Ax^*-b\|_{\infty}$&$\|(A^Tu^*-c)_+\|_{\infty} $ & $|c^Tx^*-b^Tu^*|$\\ [1mm] \hline \hline  &&&& \\ [-9mm]
\begin{tabular}{c}
   \ \\ $800,1000,1$  \end{tabular}& linprog & $$ &$ $ &$Out ~of~ memory$  &$ $& $ $ \\
   [-4mm]\cline{2-7} & lpf & $26.06$ & $1.3250e+001$ & $4.6280e-009$&$1.3358e-012 $& $ 2.2695e-008$\\ &&&& \\ [-7mm] \hline \hline &&&& \\ [-9mm]
\begin{tabular}{c}
   \ \\ $1000,1200,1$  \end{tabular}& linprog & $$ &$ $ &$Out ~of~ memory$  &$ $& $ $ \\
   [-4mm]\cline{2-7} & lpf & $45.91$ & $1.3072e+001$ & $5.8478e-009$&$1.7586e-012 $& $ -1.6425e-008$\\ &&&& \\ [-7mm] \hline \hline &&&& \\ [-9mm]
   \begin{tabular}{c}
   \ \\ $800,10000,1$ \end{tabular}& linprog & $$ &$ $ &$Out ~of~ memory$  &$ $& $ $ \\
   [-4mm]\cline{2-7} & lpf & $97.20$ & $9.0806e+000$ & $1.1331e-008$&$5.6843e-013 $& $ 8.8519e-008$\\ &&&& \\ [-7mm] \hline \hline &&&& \\ [-9mm]
   \begin{tabular}{c}
   \ \\ $2000,30000,0.1$ \end{tabular}& linprog & $162.51$ &$3.2136e+001 $ &$1.2335e-008$  &$5.79e-013 $& $5.82e-007$ \\
   [-4mm]\cline{2-7} & lpf & $29.31$ & $13.372$ & $3.9063e-009$&$6.8212e-013 $& $ -5.6716e-009$\\ &&&& \\ [-7mm] \hline \hline &&&& \\ [-9mm]
   \begin{tabular}{c}
   \ \\ $3000,4000,0.1$ \end{tabular}& linprog & $$ &$ $ &$Out ~of~ memory$  &$ $& $ $ \\
   [-4mm]\cline{2-7} & lpf & $87.35$ & $1.3993e+001$ & $8.2518e-009$&$8.5265e-013$& $ 4.9171e-008$\\ &&&& \\ [-7mm] \hline \hline &&&& \\ [-9mm]
   \begin{tabular}{c}
   \ \\ $10000,15000,0.01$  \end{tabular}& linprog  & $$ &$ $ &$Out ~of~ memory$  &$ $& $ $ \\
   [-4mm]\cline{2-7} & lpf & $$ &$ $ &$Out ~of~ memory$  &$ $& $ $ \\ &&&& \\ [-7mm] \hline \hline &&&& \\ [-9mm]
   \begin{tabular}{c}
   \ \\ $10000,150000,0.01$  \end{tabular}& linprog & $$ &$ $ &$Out ~of~ memory$  &$ $& $ $ \\
   [-4mm]\cline{2-7} & lpf & $$ &$ $ &$Out ~of~ memory$  &$ $& $ $ \\ &&&& \\ [-7mm] \hline \hline &&&& \\ [-9mm]
   \begin{tabular}{c}
   \ \\ $100,10000000,0.01$  \end{tabular}& linprog & $814.07$ &$1.1705e+005 $ &$4.0059e-007$  &$1.00e-007 $& $1.82e-001$ \\
   [-4mm]\cline{2-7} & lpf & $218.86$ & $4.3655e-001$ & $2.0084e-007$&$7.1054e-015 $& $ -2.8666e-007$ \\ [0.5mm] \hline
\end{tabular}
}
\end{center}
\end{sidewaystable}

\begin{sidewaystable}
\caption{ Comparison of \textit{\textbf{linprog}} and  \textit{\textbf{lpf}} in Linux}
\begin{center}{
\begin{tabular}{|c|c|c|c|c|c|c|}\hline &&&&\\ [-6mm]
   $m,n,d$&&$ time(sec)$&$\|{x^*}\|$& $\|Ax^*-b\|_{\infty}$&$\|(A^Tu^*-c)_+\|_{\infty} $ & $|c^Tx^*-b^Tu^*|$\\ [1mm] \hline \hline  &&&& \\ [-9mm]
\begin{tabular}{c}
  \ \\ $800,1000,1$  \end{tabular}& linprog & $45.84$ &$1.8190e+01 $ &$3.4152e-10$  &$8.01e-13 $& $4.97e-06 $ \\
   [-4mm]\cline{2-7} & lpf & $23.83$ & $1.3250e+01$ & $4.6280e-09$&$1.3358e-12 $& $ 2.2695e-08$\\ &&&& \\ [-7mm] \hline \hline &&&& \\ [-9mm]
\begin{tabular}{c}
   \ \\ $1000,1200,1$ \end{tabular}& linprog & $$ &$ $ &$Out ~of~ memory$  &$ $& $ $ \\
   [-4mm]\cline{2-7} & lpf  & $44.59$ &  $1.3072e+01$ & $5.8478e-09$&$1.7586e-12 $& $ -1.6425e-08$\\ &&&& \\ [-7mm] \hline \hline &&&& \\ [-9mm]
   \begin{tabular}{c}
   \ \\ $800,10000,1$  \end{tabular}& linprog & $$ &$ $ &$Out ~of~ memory$  &$ $& $ $ \\
   [-4mm]\cline{2-7} & lpf  & $90.41$ &  $9.0806e+00$ & $1.1331e-08$&$5.6843e-13 $& $ 8.8519e-08$\\ &&&& \\ [-7mm] \hline \hline &&&& \\ [-9mm]
   \begin{tabular}{c}
   \ \\ $2000,30000,0.1$  \end{tabular}& linprog & $161.72$ &$3.2136e+01$&$9.4724e-10$ &$5.78e-13$& $5.82e-07 $ \\
   [-4mm]\cline{2-7} & lpf  & $34.65$ &  $13.372$ & $3.9063e-09$&$6.8212e-13 $& $ -5.6716e-09$\\ &&&& \\ [-7mm] \hline \hline &&&& \\ [-9mm]
  \begin{tabular}{c}
   \ \\ $3000,4000,0.1$  \end{tabular}& linprog & $291.53$ &$2.3014e+01 $&$6.4369e-10$ &$1.02e-12 $& $9.65e-09 $ \\
   [-4mm]\cline{2-7} & lpf  & $113.95$ &  $1.3993e+01$ & $8.2518e-09$&$8.5265e-13$& $4.9171e-08$\\ &&&& \\ [-7mm] \hline \hline &&&& \\ [-9mm]
   \begin{tabular}{c}
   \ \\ $10000,15000,0.01$  \end{tabular}& linprog & $$ &$ $ &$Out ~of~ memory$  &$ $& $ $ \\
   [-4mm]\cline{2-7} & lpf  & $3243.50$ &  $1.4390e+01$ & $6.0431e-09$&$3.1264e-13 $& $ -1.2012e-07$\\ &&&& \\ [-7mm] \hline \hline &&&& \\ [-9mm]
   \begin{tabular}{c}
   \ \\ $10000,150000,0.01$ \end{tabular}& linprog & $$ &$ $ &$Out ~of~ memory$  &$ $& $ $ \\
   [-4mm]\cline{2-7} & lpf  & $4215.10$ &  $5.5538e+01$ & $1.1173e+05$&$1.1108e-01$& $ 3.8562e+05$\\ &&&& \\ [-7mm] \hline \hline &&&& \\ [-9mm]
   \begin{tabular}{c}
   \ \\ $100,10000000,0.01$  \end{tabular}& linprog & $212.72$ &$1.1705e+05 $&$4.4343e-07$ &$1.00e-07$& $1.82e-01 $ \\
   [-4mm]\cline{2-7} & lpf  & $140.01$ &  $4.3655e-01$ & $2.0084e-07$&$7.1054e-15 $& $ -2.8666e-07$ \\ [0.5mm] \hline
\end{tabular}
}
\end{center}
\end{sidewaystable}

\section{Conclusion}
~~~Two methods have been offered in this paper to solve linear programming problem. Initially, some problems are produced randomly and then they will be solved by two offered methods in the Win OS. Though the second algorithm enjoys a higher efficiency than MATLAB LINPROG according to the results, it cannot solve the large-scale problems. Then the produced problems are solved using the both methods in the Linux OS, the results show that some of the unsolved problems in the Win OS are solved in the Linux OS.\\
\\
\newpage

\end{document}